\documentclass[11pt]{amsart}
\usepackage{ae}
\usepackage[all]{xy}
\usepackage{graphicx,amsfonts,amssymb,amsmath}
\usepackage{amsmath}
\usepackage{pifont}

 \usepackage[backref]{hyperref}
\hypersetup{colorlinks,linkcolor=blue,citecolor=red,}
\usepackage[alphabetic,backrefs,lite]{amsrefs}

\usepackage[OT2,T1]{fontenc}
\DeclareSymbolFont{cyrletters}{OT2}{wncyr}{m}{n}
\DeclareMathSymbol{\sha}{\mathalpha}{cyrletters}{"58}
\input cyracc.def

 \newtheorem{thm}{Theorem}[section]

 \newtheorem{cor}[thm]{Corollary}
 \newtheorem{lem}[thm]{Lemma}
 
 \newtheorem{prop}[thm]{Proposition}
 \theoremstyle{definition}
 \newtheorem{defn}[thm]{Definition}
 \theoremstyle{remark}
 
 \theoremstyle{remark}
 \newtheorem{rem}[thm]{Remark}
 \numberwithin{equation}{subsection}
 \newcommand{\To}{\longrightarrow}

 \newcommand{\Br}{\textup{Br}}

  \newcommand{\pr}{\textup{pr}}
 
  \newcommand{\Spec}{\textup{Spec}}
 \newcommand{\BM}{Brauer\textendash Manin\ }
 \newcommand{\BMo}{Brauer\textendash Manin obstruction\ }

 \renewcommand{\P}{\mathbb{P}}

 \newcommand{\A}{\textbf{A}}
 \newcommand{\Q}{\mathbb{Q}}
 \newcommand{\R}{\mathbb{R}}
  \newcommand{\C}{\mathbb{C}}
 \newcommand{\Z}{\mathbb{Z}}

\makeatletter\let\@wraptoccontribs\wraptoccontribs\makeatother

\begin{document}

\title[]
 {Arithmetic of Ch\^atelet surface bundles revisited}

\author{Guang Hu}
\author{Yongqi Liang}

\address{Guang Hu
\newline University of Scinece and Technology of China,
\newline School of Mathematical Sciences,
\newline 96 Jinzhai Road,
\newline  230026 Hefei, Anhui, China
 }

\email{huguang@mail.ustc.edu.cn}

\address{Yongqi LIANG
\newline University of Scinece and Technology of China,
\newline School of Mathematical Sciences,
\newline 96 Jinzhai Road,
\newline  230026 Hefei, Anhui, China
 }

\email{yqliang@ustc.edu.cn}

%\date{\today.}
%\dedicatory{}

%%% ----------------------------------------------------------------------

%%% ----------------------------------------------------------------------
\maketitle

\begin{abstract}
We study arithmetic of the algebraic varieties defined over number fields by applying Lagrange interpolation to fibrations.

Assuming the finiteness of the Tate\textendash Shafarevich group of a certain elliptic curve, we show, for Ch\^atelet surface bundles over curves, that the violation of Hasse principle being accounted for by the Brauer\textendash Manin obstruction is not invariant under an arbitrary finite extension of the ground field.
\end{abstract}

%%% ----------------------------------------------------------------------

\section{Introduction}
In arithmetic geometry, people study whether diverse families of algebraic varieties defined over number fields satisfy Hasse principle. If violation of Hasse principle happens, people ask whether it is accounted for by the  Brauer\textendash Manin obstruction, cf. \cite{Skbook}, \cite{Poonenbook}, or \cite{CTSkbook} for an introduction to this topic.

The arithmetic of Ch\^atelet surfaces has been well-studied in a profound paper of J.-L. Colliot-Th\'el\`ene, J.-J. Sansuc and P. Swinnerton-Dyer \cite{chateletsurfaces0,chateletsurfaces}. They showed that Ch\^atelet surfaces may violate Hasse principle, but the violation is always accounted for by the Brauer\textendash Manin obstruction.

The study of the arithmetic of Ch\^atelet surface bundles was started by Bjorn Poonen in his paper \cite{Poonen}. He constructed a Ch\^atelet surface bundle over a curve of high genus to show that the failure of Hasse principle may not be explained by Brauer\textendash Manin obstruction even applied to \'etale covers.

We continued to study the arithmetic of such Ch\^atelet surface bundles. In the previous paper \cite{Liang-noninv} of the second author, we focused on its behaviors under the extension of the ground field. We showed, for a certain quadratic extension, that the property of weak approximation with Brauer\textendash Manin obstruction off infinity is not invariant. If the base field is a number field other than the field $\mathbb{Q}$ of rational numbers, the result is conditional on a Conjecture of Michael Stoll, cf. \cite[Conjecture 9.1]{Stoll07}.

Under the same conjecture, Han Wu extended the previous result to  arbitrary nontrivial extensions of number fields, cf. \cite[Theorem 6.2.1]{WuHan}. In the same paper, Wu also studied the behavior under extension of the ground field of the failure of Hasse principle for such Ch\^atelet surface bundles. Assuming Stoll's conjecture over the ground field $K$, he showed that this failure can be explained by Brauer\textendash Manin obstruction over $K$ while it may  not be the case after base extension to $L$, provided that the nontrivial extension $L/K$ satisfies one of the following conditions:
\begin{itemize}
\item[-] the degree $[L:K]$ is odd;
\item[-] both $K$ and $L$ admit  real places;
\end{itemize}
cf. \cite[Theorems 6.3.1 and 6.3.2]{WuHan}. Moreover, when $K=\mathbb{Q}$ and $L=\mathbb{Q}(\sqrt{-1})$ where neither condition above on $L/K$ is fulfilled, he also proved the same statement by an explicit construction \cite[\S 7.4]{WuHan}.

The main result of the present paper is the following Theorem \ref{maintheorem}. It states roughly that, assuming the finiteness of the Tate\textendash Shafarevich group of a certain elliptic curve, for an arbitrary nontrivial extension $L/K$ of number fields, the failure of Hasse principle of Ch\^atelet surface bundles over curves being accounted for by the Brauer\textendash Manin obstruction is not invariant under the extension  of the ground field. When $K=\Q$ the assumption is verified, the result turns out to be unconditional. Theorem \ref{maintheorem}  simultaneously generalises Wu's Theorem 6.3.1, Theorem 6.3.2 and \S 7.4 in \cite{WuHan}. It also says that for the \'etale Brauer\textendash Manin obstruction   similar non-invariant phenomenon happens to Ch\^atelet surface bundles as well. Our method provides a uniform way to discuss the arithmetic of Ch\^atelet surface bundles over curves, it may have possible applications to other fibrations.

Apart from Ch\^atelet surface bundles, we should also mention that such non-invariant phenomenon for certain surfaces has also been studied by Wu in his subsequent paper \cite{WuHan2}.

Note that the study of invariance of arithmetic properties under  extension of the ground field arises naturally. In \cite[Corollary 1.7]{BN2021}, for Kummer varieties, Francesca Balestrieri and Rachel Newton proved that the existence of Brauer\textendash Manin obstruction to Hasse principle may be preserved under  infinitely many possible extensions of the ground field.
A similar statement also appeared in Manar Riman's preprint \cite{Riman} for del Pezzo surfaces of degree $4$, see also the recent paper \cite[Lemma 7.2]{CreutzViray21} of Brendan Creutz and Bianca Viray. In \cite{RiveraViray21}, Carlos Rivera and Bianca Viray study similar questions for cubic surfaces.
In this paper, we deal with Ch\^atelet surface bundles and our result goes completely in the other direction.

\begin{thm}\label{maintheorem}
Let $K$ be a number field. We assume the existence of an elliptic curve over $K$  of Mordell\textendash Weil rank zero such that the subgroup of divisible elements of its Tate\textendash Shafarevich group is trivial. (For example, this is the case for $K=\Q$.)

Then for any nontrivial  finite extension  $L/K$  there exists a Ch\^atelet surface bundle $X$ defined over $K$ violating the Hasse principle over both $K$ and $L$ such that
\begin{itemize}
\item[-] the \BMo is the only obstruction to Hasse principle for $X_K$, i.e. $$\varnothing=X(K)=X(\textup{\textbf{A}}_K)^{\textup{et},\Br}=X(\textup{\textbf{A}}_K)^{\Br(X)}\subsetneq X(\textup{\textbf{A}}_K),$$
\item[-] while the \'etale \BMo is not the only obstruction to Hasse principle for $X_L$, i.e. $$\varnothing=X_L(L)\subsetneq X_L({\textup{\textbf{A}}}_L)^{\textup{et},\Br}\subset X_L({\textup{\textbf{A}}}_L)^{\Br(X_L)}\subset X_L({\textup{\textbf{A}}}_L).$$
\end{itemize}
\end{thm}

We remark that Stoll's conjecture for an elliptic curve is equivalent to the condition that the subgroup of divisible elements of its Tate\textendash Shafarevich group is trivial cf. \cite[Corollary 6.2]{Stoll07}, which is a consequence of the famous conjecture on
the finiteness of its Tate\textendash Shafarevich group. In the proofs of previous results in \cite{Liang-noninv} and \cite{WuHan}, the assumption of Stoll's conjecture on higher genus curves were required, while here we succeed  to replace the assumption by the one on elliptic curves. In consequence, Theorem \ref{maintheorem} turns out to be unconditional for arbitrary finite extensions $L/\Q$. This recovers all unconditional examples in \cite{Liang-noninv} and \cite{WuHan}, where only very particular extensions $L$ of $\Q$ were considered.

The construction of the bundle in Theorem \ref{maintheorem} dates back to Poonen.
He began with a Ch\^atelet surface bundle over $\mathbb{P}^1$ with a particular key fiber carefully chosen over a rational point. Later in Wu's consideration, two key fibers of the bundle over $\mathbb{P}^1$ were taken care of, which were still lying over rational points. Our main ingredient is to apply Lagrange interpolation to obtain algebraic variety bundles over $\mathbb{P}^1$, so that arbitrary finite number of key fibers can be controlled at the same time. In Wu's argument \cite{WuHan}, he has to require the  extension $L/K$ to be of odd degree. Because the relevant Brauer group is of exponent $2$ but the key fibers are defined over the ground field $K$, the Brauer\textendash Manin pairing will vanish after an even degree extension and Wu's argument fails without the parity assumption. While in our treatment, key fibers can be chosen to be defined over $L$ (or any other finite extension of $K$), i.e.  more than only two geometric fibers are considered at the same time. Exactly in such a case, Wu succeeded to obtain his very particular single example for the extension $\Q(\sqrt{-1})/\Q$ in \cite[\S 7.4]{WuHan}.

The paper is organised as follows. In \S \ref{preliminaries}, we set up the background by recalling some basic facts about Ch\^atelet surface bundles. For completeness of the paper, we also recall local calculations on Ch\^atelet surfaces essentially done by Wu in \cite{WuHan}. In \S \ref{proof}, after explaining Poonen's construction, we discuss Lagrange interpolation applied to Ch\^atelet surface bundles and complete the proof of our main theorem. Finally in \S \ref{further}, we discuss the opposite situation of our main theorem.

\section{Preliminaries}\label{preliminaries}
\subsection{Notation}
We briefly recall some notation. In this paper, $L/K$ is a nontrivial extension of  number fields. For $F=K$ or $L$, the set of places (respectively archimedean places, non-archimedean places, $2$-adic places) of $F$ is denoted by $\Omega_F$ (respectively $\Omega_F^\infty$,$\Omega_F^\textup{f}$, $\Omega^2_F$).
%For $v\in\Omega_F$, $F_v$ is the completion of $F$ at $v$.
The ring of ad\`eles is denoted by $\textup{\textbf{A}}_F$.

We denote by $X$  an algebraic variety (separated scheme of finite type) over $K$.  We denote by $X_L$ the base change $X\times_{\Spec(K)}\Spec(L)$ of $X$ to $L$. $\Br(X)=\textup{H}^2_{\scriptsize{\textup{\'et}}}(X,\mathbb{G}_\textup{m})$ denotes the cohomological Brauer group of $X$.

\subsection{Ch\^atelet surface bundles}\label{ChSurBun-section}
In order to fix notation, we present a general geometric setup of Ch\^atelet surface bundles over $\P^1$.

Let $K$ be a field. Start with an affine surface $\mathcal{S}^o$ defined by $$y^2-az^2=P(x)$$ in $\mathbb{A}^3$, where $a\in K$ is a nonzero constant and $P(x)\in K[x]$ is a polynomial of degree $4$. It has a projective model $\mathcal{S}$ obtained by gluing
\begin{itemize}
\item[-] the surface defined by $Y^2-aZ^2=W^2P(x)$ in $(\P^1\setminus\{\infty\})\times\P^2$ with coordinates $(x,(Y:Z:W))$ and
\item[-] the surface defined by $Y^2-aZ^2=W'^2P^*(x')$ in $(\P^1\setminus\{0\})\times\P^2$ with coordinates $(x',(Y:Z:W'))$
\end{itemize}
via identifications of variables $x=1/x'$ and $W=x'^2W'$, where $P^*(x')=x'^4P(1/x')\in K[x']$ denotes the reciprocal polynomial of $P(x)$.
According to Jacobian criterion, it is smooth if (and only if) $P(x)$ is separable (then $P^*(x')$ is also separable and is of degree $3$ or $4$). We call this smooth projective variety $\mathcal{S}$ a \emph{Ch\^atelet surface}. It can be viewed as a conic bundle over $\P^1$ via the variable $x$ (or $x'$ around infinity).

We consider a family of Ch\^atelet surfaces parameterised over $\mathbb{A}^1$ by $t$
$$y^2-a_tz^2=P_t(x)$$
where $a_t\in K[t]$ is nonzero and $P_t(x)\in K[t,x]$ is a polynomial of degree $4$ in the variable $x$.
Replacing $a$ by $a_t$, $P(x)$ by $P_t(x)$, and $P^*(x')$ by $P^*_t(x')$ respectively in the projective model, we obtain a hypersurface in $\mathbb{A}^1\times\P^1\times\P^2$ and we call it a \emph{Ch\^atelet surface bundle} over $\mathbb{A}^1$. For almost all algebraic values $\theta$ of $t$ (strictly speaking,  closed points $\theta$ of $\mathbb{A}^1$), the fiber $$V_\theta: \quad y^2-a_\theta z^2=P_\theta(x)$$ is a Ch\^atelet surface defined over the residue field $K(\theta)$.

We want to extend such a Ch\^atelet surface bundle over $\mathbb{A}^1$ to a bundle $V$ over $\P^1$. In order to have an explicit projective model completely determined by the polynomial $P_t(x)$, we need two further assumptions:
\begin{itemize}
\item[-] the polynomial $a_t\in K[t]$ is a nonzero constant $a\in K$,
\item[-] the polynomial $P_t(x)\in K[t,x]$ is of \emph{even} degree $d$ in the variable $t$.
\end{itemize}

The variety $V$ will contain a Zariski dense open subset $V^o$ which is a hypersurface  of $\mathbb{A}^1\times\mathbb{A}^1\times\P^2$ defined by the equation
$$Y^2-aZ^2=W^2P_{t}(x)$$ with coordinates $(t,x,(Y:Z:W))$. It is a conic bundle over $\mathbb{A}^1\times\mathbb{A}^1$ via the variables $(t,x)$. We are going to extend it to a conic bundle $V\To\P^1\times\P^1$, and the Ch\^atelet surface bundle structure can be seen by applying a further projection $\pr_1$ to the first factor.

Let $\P^1\times\P^1$ be covered by affine open subsets
\begin{equation*}
    \begin{array}{rl}
        \mathcal{U}_{\infty\infty} & =(\P^1\setminus\{\infty\})\times(\P^1\setminus\{\infty\}),\\
        \mathcal{U}_{\infty0} & =(\P^1\setminus\{\infty\})\times(\P^1\setminus\{0\}),\\
        \mathcal{U}_{0\infty} & =(\P^1\setminus\{0\})\times(\P^1\setminus\{\infty\}),\\
    \end{array}
\end{equation*}
and
\begin{equation*}
    \begin{array}{rl}
        \mathcal{U}_{00}  & =(\P^1\setminus\{0\})\times(\P^1\setminus\{0\}).
    \end{array}
\end{equation*}
Consider four similar equations:
\begin{equation*}
\begin{array}{lll}
\mbox{over }  \mathcal{U}_{\infty\infty}&: & Y^2-aZ^2=W_{\infty\infty}^2P^{\infty\infty}_{t}(x);\\
\mbox{over }  \mathcal{U}_{\infty0}&: & Y^2-aZ^2=W_{\infty0}^2P^{\infty0}_{t}(x');\\
\mbox{over }  \mathcal{U}_{0\infty}&: & Y^2-aZ^2=W_{0\infty}^2P^{0\infty}_{t'}(x);\\
\mbox{over }  \mathcal{U}_{00}&: & Y^2-aZ^2=W_{00}^2P^{00}_{t'}(x');
\end{array}
\end{equation*}
where the four polynomials on the right hand side are obtained from $P_t(x)$ as follows.
Let us denote by $$\tilde{P}_{t,t'}(x,x')=x'^4t'^dP_{t/t'}(x/x')\in K[t,t',x,x']$$ the homogenisation
%(separably in $t$ and in $x$, but not  in the total degree of $t$ and  $x$)
of $P_t(x)$. We set
\begin{itemize}
\item[-] $P^{\infty\infty}_{t}(x)=\tilde{P}_{t,1}(x,1)=P_t(x)$,
\item[-] $P^{\infty0}_{t}(x')=\tilde{P}_{t,1}(1,x')$ \textemdash~ the reciprocal polynomial of $P_t(x)$ w.r.t. $x$,
\item[-] $P^{0\infty}_{t'}(x)=\tilde{P}_{1,t'}(x,1)$ \textemdash~ the reciprocal polynomial of $P_t(x)$ w.r.t. $t$,
\item[-] $P^{00}_{t'}(x')=\tilde{P}_{1,t'}(1,x')$ \textemdash~ the reciprocal polynomial of $P_t(x)$ w.r.t. $t$ and $x$.
\end{itemize}
Identifying $W$ with $W_{\infty\infty}$, the first piece is exactly $V^o$.
These four pieces can be glued together to obtain  $V\To\P^1\times\P^1$ via compatible identifications of variables
$$t'=t^{-1},~ x'=x^{-1},$$
and
\begin{equation*}
    \begin{array}{rclrclrcl}
       W_{\infty\infty}&=&W_{\infty0}x'^2,&W_{\infty\infty}&=&W_{0\infty}t'^{d/2},&W_{\infty\infty}&=&W_{00}x'^2t'^{d/2},\\
       W_{\infty0}&=&W_{0\infty}x^2t'^{d/2},& W_{\infty0}&=&W_{00}t'^{d/2},& W_{0\infty}&=&W_{00}x'^2.$$
    \end{array}
\end{equation*}
We remark that the parity assumption on $d$ is applied here.

To see the Ch\^atelet surface bundle structure, we take a closed point $\theta$ of $\P^1\setminus\{\infty\}$.  By specialising the first two equations at $t=\theta$, the fiber $V_\theta$ is thus given by
 $$Y^2-aZ^2=W_{\infty\infty}^2P^{\infty\infty}_\theta(x)=W^2P_\theta(x)$$
in $(\mathbb{P}^1\setminus\{\infty\})\times\mathbb{P}^2$ with coordinates $(x,(Y:Z:W_{\infty\infty}))$ together with
 $$Y^2-aZ^2=W_{\infty0}^2P^{\infty0}_\theta(x')$$
in $(\mathbb{P}^1\setminus\{0\})\times\mathbb{P}^2$ with coordinates $(x',(Y:Z:W_{\infty0}))$. Since $P^{\infty0}_\theta(x')=P_\theta^*(x')$, this is exactly the Ch\^atelet surface bundle over $\mathbb{A}^1$ described previously when $a_t=a$ is a constant. For simplicity, we usually only write down the corresponding affine equation  $$y^2-az^2=P_\theta(x),$$ where  $P_\theta(x)\in K(\theta)[x]$. In a uniform way, for $\theta=\infty\in\P^1$, roughly speaking the fiber $V_\infty$ is  defined by  $$y^2-az^2=P_\infty(x),$$ where $P_\infty(x)$ is obtained by substituting $t'=0$ (\textup{i.e.} $t=\infty$) in the polynomial $t'^dP_{1/t'}(x)\in K[t',x]$. Precisely, one can obtain the defining equations of $V_\infty$ by substituting $t'=0$  in the equations over $\mathcal{U}_{0\infty}$ and $\mathcal{U}_{00}$. For almost all values $\theta$ of $t$, the polynomial $P_\theta(x)$ is separable of degree $4$ signifying that $V_\theta$ is a Ch\^atelet surface.

To summarize, under the assumptions that $a_t=a$ is a constant and the degree $d$  is even, the Ch\^atelet surface bundle $V\To\P^1$ is determined by the polynomial $P_t(x)\in K[t,x]$. Therefore, we only need to remember its affine equation $$y^2-az^2=P_t(x).$$

The \emph{degeneracy locus} $\textup{De}(V)\subset\P^1\times\P^1$ of the conic bundle $V\To\P^1\times\P^1$ is a projective curve defined by the equation $\tilde{P}_{t,t'}(x,x')=0.$
The following lemma is a consequence of the Jacobian criterion.

\begin{lem}\label{totalsmoothnesslemma}
When $a_t=a$ is a constant and the degree $d$ is even,  the projective variety $V$ is smooth over $K$ if and only if $\textup{De}(V)$ is a smooth curve.
\end{lem}

\begin{defn}\label{def-adm}
We say that the polynomial $P_t(x)\in K[t,x]$ is \emph{admissible} if the corresponding degeneracy locus $\textup{De}(V)\subset \P^1\times\P^1$ is a smooth curve.
\end{defn}

The following lemma also follows from the Jacobian criterion.
\begin{lem}\label{smoothfiberlemma}
For each closed point $\theta\in\P^1$, the fiber $V_\theta$ is smooth over the residue field $K(\theta)$ if and only if $\theta$ does not lie in the branch locus $\mathfrak{B}\subset\P^1$ of $$\pr_1:\textup{De}(V)\subset\P^1\times\P^1\To\P^1.$$
This is also equivalent to the condition that $P_\theta(x)\in K(\theta)[x]$ and $P_\theta^*(x')=x'^4P_\theta(1/x')\in K(\theta)[x']$ are both separable polynomials.
\end{lem}
We remark that the separability of $P^*_\theta$ follows  automatically from the separability of $P_\theta$ if $\deg(P_\theta)=4$.

\begin{defn}\label{degeneracylocusofChsurfbun}
According to the previous lemma, we call the branch locus $\mathfrak{B}\subset\P^1$ of $$\pr_1:\textup{De}(V)\subset\P^1\times\P^1\To\P^1$$ \emph{the degeneracy locus} of the Ch\^atelet surface bundle $V\To\P^1$.
\end{defn}

\subsection{Hyperelliptic curves}\label{hyperelliptic}
Let $K$ be a field of characteristic $0$.
Consider an affine plane  curve defined by the equation $$s^2=H(t)=a_nt^n+a_{n-1}t^{n-1}+\cdots+a_1t+a_0$$ with variables $s$ and $t$,
where $H\in K[t]$ is a  polynomial of  degree $n=2g+1>4$ or $n=2g+2>4$. It is a ramified double cover of $\mathbb{A}^1$ via the projection $\gamma:(s,t)\mapsto t$. Let $H^\dag\in K[t']$ denote the polynomial $H^\dag(t')=t'^{2g+2}H(1/t')$. Consider the affine plane curve defined by $$s'^2=H^\dag(t'),$$
which is also a ramified double cover of $\mathbb{A}^1$ via $\gamma':(s',t')\mapsto t'$.
We can glue these two affine curves on the common open subset given by $t\neq0$
and $t'\neq0$ via the identification of variables $t'=1/t$ and $s'=t^{-g-1}s$. When $H(t)$ is a separable polynomial, the projective curve $C$ obtained in this way is smooth and of genus $g\geq2$, we call it a \emph{hyperelliptic curve}. We often refer to  $s^2=H(t)$ as its defining equation.
We also obtain a ramified double cover $\gamma:C\To\P^1$. Ramifications of $\gamma$ appear over roots of $H$, and in addition over the point of infinity when $n$ is odd.

\subsection{Local calculations on Ch\^atelet surfaces}\label{local-ChSur}
According to \cite[Theorem 8.11]{chateletsurfaces} the existence of global rational points on a Ch\^atelet surface is determined by its Brauer\textendash Manin set, which can be computed explicitly.
A very general local calculations on Ch\^atelet surfaces is due to Han Wu \cite[\S 3]{WuHan}, of which only a small part is relevant for this paper. For completeness, we recall the crucial ones.

We recall basic properties about the Hilbert symbol. Let $F$ be a completion of a number field at a place $v$. For $\alpha,\beta\in F^*$, the \emph{Hilbert symbol} $(\alpha,\beta)_v\in\{\pm1\}$ is defined to be $1$ if the equation $x^2-\alpha y^2-\beta z^2=0$ has a nonzero solution in $F$, otherwise it is defined to be $-1$. We know that the Hilbert symbol is $\mathbb{F}_2$-bilinear and symmetric, cf. \cite[Chapitre XIV Proposition 7]{SerreCorpsLoc}.

\begin{lem}\label{lem5}
Suppose that $v$ is an odd place. Then  $(\alpha,\beta\pm\beta')_v=(\alpha,\beta)_v$ if $v(\beta)<v(\beta')<+\infty$.
\end{lem}
\begin{proof}
As $(\alpha,\beta\pm\beta')_v=(\alpha,\beta(1\pm\beta'\beta^{-1}))_v=(\alpha,\beta)_v\cdot(\alpha,1+\beta'\beta^{-1})_v$, it remains to show that $(\alpha,1+\beta'\beta^{-1})_v=1$. But it  follows from Hensel's lemma that the $v$-adic integer $1+\beta'\beta^{-1}$ is a square.
\end{proof}

\begin{cor}\label{cor7}
Suppose that $v$ is an odd place, then the map $\beta\mapsto (\alpha,\beta)_v$ is locally constant.
\end{cor}

\begin{lem}\label{lem6}
Suppose that $v$ is an odd place. If $(\alpha,\beta)_v=-1$ then at least one of $v(\alpha)$ and $v(\beta)$ is an odd integer.

Conversely, if $v(\alpha)$ is odd, then there exists $\beta\in F$ such that $(\alpha,\beta)_v=-1$, moreover one can require that $v(\beta)=0$.
\end{lem}
\begin{proof}
If both $v(\alpha)$ and $v(\beta)$ are even, it reduces to the case where they are all $0$. But Hensel's lemma and Chevalley\textendash Warning's theorem applied to quadratic forms in three variables imply that $(\alpha,\beta)_v=1$ leading to a contradiction.

Conversely, as $\alpha$ can never be a square in $F$, the existence of $\beta$ follows from \cite[Chapitre XIV Proposition 7 vi)]{SerreCorpsLoc}.
If $v(\beta)$ is odd, then $v(-\alpha\beta)$ is even and $(\alpha,-\alpha\beta)_v=(\alpha,-\alpha)_v\cdot(\alpha,\beta)_v=-1$ by \cite[Chapitre XIV Proposition 7 iv)]{SerreCorpsLoc}. Hence we may require that $v(\beta)$ is even, and further more $v(\beta)=0$ by multiplying a suitable even power of a uniformizer of $F$.
\end{proof}

For the convenience of exposition, we make the following definition in this paper.
\begin{defn}\label{2Rlocal}
For a number field $K$, we say that an element $a\in K$ is an \emph{$2\mathbb{R}$ local square} if it is a square in the completion $K_v$ for all even places and for all real places.
\end{defn}
It is clear that if $a\in K$ is an $2\mathbb{R}$ local square, then it is also an $2\mathbb{R}$ local square viewed as an element in any finite extension $L$ of $K$.

For the convenience of the reader, we make the following proposition of Han Wu adapted to our proof.

\begin{prop}[compare to {\cite[Proposition 3.1.1]{WuHan}}]\label{WuHan1}
Let $K$ be a number field and $a\in\mathcal{O}_K$ be a non-zero $2\mathbb{R}$ local square (cf. Definition \ref{2Rlocal}).

Then for any non-archimedean place $v_0$ of $K$ such that $v_0(a)$ is odd, there exists a Ch\^atelet surface $\mathcal{S}$ defined over $K$ by an equation
$$y^2-az^2=Ax^4+Cx^2+E$$ with $E\neq0$, such that $\mathcal{S}(K_{v_0})=\varnothing$ but $\mathcal{S}(K_{v})\neq\varnothing$ for all $v\neq v_0$.

In particular, for any finite extension $L$ of $K$ on which $v_0$ splits completely, the surface $\mathcal{S}$ does not possess any $L$-rational point.
\end{prop}

\begin{proof}
When $v_0$ splits completely in $L$, the field $L_{w_0}$ equals to $K_{v_0}$ for any place $w_0$ of $L$ lying over $v_0$. The last statement follows from $\mathcal{S}(K_{v_0})=\varnothing$.

It remains to prove the main statement.
The assumption implies that $v_0$ is an odd place.
Define $S_1=\{v\in\Omega_K^\text{f}; v\nmid2\mbox{ and }v(a)\neq0~\}$, then $v_0\in S_1$. By weak approximation of $K$, we may choose a constant $b\in K$ such that
\begin{itemize}
\item[-] $(a,b)_{v_0}=-1$   (this is a nonempty open condition on $K_{v_0}$ by Lemma \ref{lem6} and Corollary \ref{cor7}),
\item[-] $(a,b)_v=1$ for all $v\in S_1\setminus\{v_0\}$.
\end{itemize}
Define $S_2=\{v\in\Omega_K^\text{f};~v\nmid2\mbox{ and }v(b)\mbox{ is odd}\}$. Again by weak approximation of $K$, we may choose a constant $c\in K$ such that
\begin{itemize}
\item[-] $c\in K_v^{*2}$ for all $v\in S_1$,
\item[-] $v(c)$ is odd for all $v\in S_2\setminus S_1$.
\end{itemize}
Since the properties listed above are stable under multiplication by a square in $K^*$, we may assume that $b,c\in\mathcal{O}_K$.

Now consider the affine surface $\mathcal{S}^o$ defined by $$y^2-az^2=b(x^4-ac),$$ whose projective model is a Ch\^atelet surface $\mathcal{S}$ of the required form with  $A=b$, $B=C=D=0$, and $E=-abc\neq0$.

We are going to verify that $\mathcal{S}$ possesses rational points everywhere locally except at $v_0$. If not otherwise specified, every step in the forthcoming local calculations follows either from Lemma \ref{lem5} or the $\mathbb{F}_2$-bilinearity of the Hilbert symbol.
\begin{itemize}
\item[-] For $v\in\Omega_K^{\infty}\cup\Omega_K^2$, as $a$ is a square in $K_v$ the surface $\mathcal{S}$ is rational over $K_v$.
\item[-] For $v\in S_1\setminus\{v_0\}$, we have $(a,b)_v=1$ and $v(ac)\geq0$. Then for $x_0\in K_v$ such that $v(x_0)<0$, we find that $(a,b(x_0^4-ac))_v=(a,b)_v\cdot(a,x_0^4-ac)_v=(a,x_0^4-ac)_v=(a,x_0^4)_v=1$. This implies that $\mathcal{S}$ possesses a $K_v$-rational point with coordinate $x=x_0$.
\item[-] For $v\in S_2\setminus S_1$, we have $v(a)=0$ and $v(bc)$ is even. Then $(a,b(0^4-ac))_v=(a,-abc)_v=1$ by Lemma \ref{lem6}. This implies that $\mathcal{S}$ possesses a $K_v$-rational point with coordinate $x=0$.
\item[-] For $v\in\Omega_K^\text{f}\setminus(S_1\cup S_2\cup\Omega^2_K)$, we have $v(ac)\geq0$, $v(a)=0$ and $v(b)$ is even, therefore $(a,b)_v=1$. Then for $x_0\in K_v$ such that $v(x_0)<0$, we find that $(a,b(x_0^4-ac))_v=(a,b)_v\cdot(a,x_0^4-ac)_v=(a,x_0^4-ac)_v=(a,x_0^4)_v=1$. This implies that $\mathcal{S}$ possesses a $K_v$-rational point with coordinate $x=x_0$.
\item[-] For $v=v_0$, concerning $x_0\in K_{v_0}$ we have $(a,b(x_0^4-ac))_{v_0}=(a,b)_{v_0}\cdot(a,x_0^4-ac)_{v_0}=-(a,x_0^4-ac)_{v_0}$.
Since $c$ is a square in $K_{v_0}$, the odd number $v_0(ac)$ can never equal to $v(x_0^4)$. If $4v_0(x_0)<v_0(ac)$, then the last equality continues to $=-(a,x_0^4)_{v_0}=-1$. If $4v_0(x_0)>v_0(ac)$, then the last equality continues to $=-(a,-ac)_{v_0}=-(a,-a)_{v_0}=-1$, where the last equality is well-known \cite[Chpitre XIV Proposition 7 iv)]{SerreCorpsLoc}.
These two cases together imply that the affine surface $\mathcal{S}^o$ does not possess any $K_{v_0}$-rational points. Since $\mathcal{S}$ is smooth and contains $\mathcal{S}^o$ as a Zariski dense open subset, the implicit function theorem implies that $\mathcal{S}$ does not possess any $K_{v_0}$-rational points either.
\end{itemize}
In summary,  $\mathcal{S}(K_{v_0})=\varnothing$ but $\mathcal{S}(K_{v})\neq\varnothing$ for all $v\neq v_0$.
\end{proof}

Over any number field, the existence of a Ch\^atelet surface violating the Hasse principle was proved by Poonen in \cite[Proposition 5.1]{Poonen09}.
In \cite{WuHan}, Han Wu did some explicit local calculations to study the invariance of the arithmetic of certain Ch\^atelet surfaces under a base extension $L/K$. Though it is not stated in the following way,  the proof of the following proposition is essentially due to Wu. For the sake of our proof,  his argument is adapted as follows.

\begin{prop}[compare to {\cite[Proposition 3.3.1]{WuHan}}]\label{WuHan2}
Let $L/K$ be an extension of number fields. Let $a\in\mathcal{O}_K$ be a non-zero $2\mathbb{R}$ local square (cf. Definition \ref{2Rlocal}).
If there exists a non-archimedean place $w_0$ of $L$ such that $w_0(a)$ is odd, then there exists a Ch\^atelet surface $\mathcal{S}$ defined over $L$ by an equation
$$y^2-az^2=Ax^4+Cx^2+E$$ with $E\neq0$, violating  the Hasse principle over $L$.
\end{prop}

\begin{proof}
Since the nonzero integer $a\in\mathcal{O}_K\subset\mathcal{O}_L$ is an $2\mathbb{R}$ local square over $L$,
the place $w_0$ must be an odd place.

Define $S_1=\{w\in\Omega_L^\text{f};~w\nmid 2 \mbox{ and }w(a)\neq0\}$. It is a finite set containing the place $w_0$. By strong approximation of $L$ without control on archimedean places, we may choose a constant $b\in\mathcal{O}_{L}$ such that
\begin{itemize}
\item[-] $w_0(b)=0$ and $(a,b)_{w_0}=-1$ (this is a nonempty open condition on $L_{w_0}$ by Lemma \ref{lem6} and Corollary \ref{cor7}),
\item[-] $w(b)=0$ and $(a,b)_w=1$ for all $w\in S_1\setminus \{w_0\}$.
\end{itemize}
Define $S_2=\{w\in\Omega_L^\text{f};~w\nmid 2 \mbox{ and }w(b)\neq0\}$. Then $S_1\cap S_2=\varnothing$ by the choice of $b$. By the Chebotarev density theorem, we choose two different non-archimedean places $w_1$ and $w_2$ outside $S_1\cup S_2\cup\Omega^2$ such that they split in the quadratic extension $L(\sqrt{a})/L$. Again by strong approximation of $L$ without control on archimedean places, we may choose a constant $c\in\mathcal{O}_{L}$ such that
\begin{itemize}
\item[-] $w(bc+1)\geq3$ for all $w\in S_1$,
\item[-] $(a,c)_w=1$ for all $w\in S_2$,
\item[-] $w_1(c)=1$,
\item[-] $w_2(bc+1)=1$.
\end{itemize}
We deduces that
\begin{itemize}
\item[-] $a$ is a unit outside $S_1\cup\Omega^2_L\cup\Omega_L^\infty$,
\item[-] $b$ is a unit outside $S_2\cup\Omega^2_L\cup\Omega_L^\infty$,
\item[-] $c$ is a unit at places in $S_1$.
%\item[-] $a$ divides $bc+1$ in $\mathcal{O}_{K,\Omega^2}$  ??useful??.
\end{itemize}

Now consider the affine surface $\mathcal{S}^o$ defined by $$y^2-az^2=(x^2-c)(bx^2-bc-1).$$
The polynomials $x^2-c$ and $bx^2-bc-1$ are irreducible over $L$ according to Eisenstein's criterion applied to $w_1$ and $w_2$ respectively.
Since they have no common root, the polynomial $(x^2-c)(bx^2-bc-1)$ is separable. The associated projective model $\mathcal{S}$ is a Ch\^atelet surface of the required form with constants $B=D=0$, $A=b$, $C=-2bc-1$, and $E=c(bc+1)\neq0$.
Moreover, the place $w_1$ is ramified in the quadratic extension $L[x]/(x^2-c)$ of $L$ and the place $w_2$ is ramified in the quadratic extension $L[x]/(bx^2-bc-1)$. Since neither of these two extension equals to the quadratic extension $L(\sqrt{a})$ where $w_1$ and $w_2$ splits, they are linearly disjoint from $L(\sqrt{a})$ over $L$. Skorobogatov's result \cite[Propositions 7.1.1 \& 7.1.2]{Skbook} states that the class of the quaternion algebra $A=(a,x^2-c)=(a,bx^2-bc-1)\in\Br(L(\mathcal{S}))$ actually lies in $\Br(\mathcal{S})$ and it generates $\Br(\mathcal{S})/\Br(L)\simeq\mathbb{Z}/2\mathbb{Z}$.

We are ready to verify that $\mathcal{S}$ violates the Hasse principle. If not otherwise specified, every step in the forthcoming local calculations follows either from Lemma \ref{lem5} or the $\mathbb{F}_2$-bilinearity of the Hilbert symbol.

First, we verify that $\mathcal{S}$ has local rational points everywhere.
\begin{itemize}
\item[-] For $w\in\Omega_L^{\infty}\cup\Omega_L^2$, as $a$ is a square in $L_w$ the surface $\mathcal{S}$ is rational over $L_w$.
\item[-] For $w=w_0$, let $\pi_0\in\mathcal{O}_{L_{w_0}}$ be a uniformizer. We have $w_0(b)=w_0(c)=0$, $w_0(\pi_0^2)=2$, $w_0(b\pi_0^2)=2$, and $w_0(bc+1)\geq3$. We find that $(a,(\pi_0^2-c)(b\pi_0^2-bc-1))_{w_0}=(a,\pi_0^2-c)_{w_0}\cdot(a,b\pi_0^2-bc-1)_{w_0}=(a,-c)_{w_0}\cdot(a,b\pi_0^2)_{w_0}=(a,-bc)_{w_0}=(a,1-(bc+1))_{w_0}=(a,1)_{w_0}=1$.
        This implies that $\mathcal{S}$ possesses a $L_{w_0}$-rational point with coordinate $x=\pi_0$.
\item[-] For $w\in S_2$, we have $w(b)>0$ and $w(c)\geq0$.  We find that $(a,(0^2-c)(b\cdot0^2-bc-1))_w=(a,c(bc+1))_w=(a,c)_w\cdot(a,1+bc)_w=(a,c)_w\cdot(a,1)_w=1$. This implies that $\mathcal{S}$ possesses a $L_{w}$-rational point with coordinate $x=0$.
\item[-] For $w\in \Omega^\text{f}_L\setminus(\{w_0\}\cup S_2\cup\Omega^2_L)$, we choose an element $x_0\in K_w$ such that $w(x_0)<0$. We have $w(b)=0$, $w(c)\geq0$, $w(bc+1)\geq0$.  We find that $(a,(x_0^2-c)(bx_0^2-bc-1))_w=(a,x_0^2-c)_w\cdot(a,bx_0^2-bc-1)_w=(a,x_0^2)_w\cdot(a,bx_0^2)_w=(a,b)_w=1$, where the last equality follows from the choice of $b$ for $w\in S_1\setminus\{w_0\}$ and otherwise from Lemma \ref{lem6}. This implies that $\mathcal{S}$ possesses a $L_{w}$-rational point with coordinate $x=x_0$.
\end{itemize}

Second, we compute the local invariant $\textup{inv}_w(A(-))\in \Q/\mathbb{Z}$ of the evaluation  of $A\in\Br(\mathcal{S})$ on $\mathcal{S}(K_w)$. Since the evaluation is locally constant, we only need to consider local points in the Zariski dense affine open subset $\mathcal{S}^o$.
\begin{itemize}
\item[-] For $w\in\Omega_L^{\infty}\cup\Omega_L^2$, the element $a$ is a square in $L_w$. Then the class of the quaternion $A$ is trivial in $\Br(\mathcal{S}_{L_w})$, its evaluation is identically $0$.
\item[-] For $w=w_0$, we have $w_0(b)=w_0(c)=0$ and $w_0(bc+1)\geq3$.  If $w_0(x)\leq0$,  we find that $(a,bx^2-bc-1)_{w_0}=(a,bx^2)_{w_0}=(a,b)_{w_0}=-1$. If $w_0(x)>0$,  we find that $(a,x^2-c)_{w_0}=(a,-c)_{w_0}=-(a,b)_{w_0}\cdot(a,-c)_{w_0}=-(a,-bc)_{w_0}=-(a,-bc+(bc+1))_{w_0}=-(a,1)_{w_0}=-1$. In both cases $\textup{inv}_{w_0}(A(-))=1/2\in\Q/\Z$.
\item[-] For $w\in S_1\setminus\{w_0\}$, we have $w(b)=w(c)=0$ and $w(bc+1)\geq3$. If $w(x)=0$, we find that $(a,bx^2-bc-1)_w=(a,bx^2)_w=(a,b)_w=1$. If $w(x)<0$, we find that $(a,x^2-c)_w=(a,x^2)_w=1$. If $w(x)>0$, we find that $(a,x^2-c)_w=(a,-c)_w=(a,b)_w\cdot(a,-c)_w=(a,-bc)_w=(a,-bc+(bc+1))_w=(a,1)_w=1$. In all cases $\textup{inv}_{w}(A(-))=0\in\Q/\Z$.
\item[-] For $w\in \Omega^\text{f}_L\setminus(S_1\cup\Omega^2_L)$,  we have $w(a)=0$, $w(b)\geq0$, and $w(c)\geq0$. If $w(x^2-c)$ is even, then $(a,x^2-c)_w=1$ by Lemma \ref{lem6}. If $w(x^2-c)$ is odd, we must have $w(x)\geq0$ and $w(x^2-c)>0$. It follows that $w(bx^2-bc-1)=w(b(x^2-c)-1)=0$ and then $(a,bx^2-bc-1)_w=1$ again by Lemma \ref{lem6}. In both cases, we find that $\textup{inv}_{w}(A(-))=0\in\Q/\Z$.
\end{itemize}
In summary, for any family of local rational points $(P_w)_{w\in\Omega_L}\in\mathcal{S}(\textup{\textbf{A}}_L)$ we have always
$$\sum_{w\in\Omega_L}\textup{inv}_w(A(P_w))=\frac{1}{2}\in\Q/\Z.$$
Whence, the surface $\mathcal{S}$ does not possess any $L$-rational point, it violates the Hasse principle.
\end{proof}

\section{Proof of the main result}\label{proof}

The aim of this section is to give a complete proof of our main result Theorem \ref{maintheorem}. It divides into several steps.

\subsection{Poonen's construction}\label{poonenconstruction}
We recall Poonen's construction of certain Ch\^atelet surface bundles over  smooth projective curves. In his paper \cite{Poonen}, B. Poonen showed that the violation of Hasse principle for such bundles may not be accounted for by the Brauer\textendash Manin obstruction even applied to \'etale covers.

Let $K$ be a number field. Consider a Ch\^atelet surface bundle $$f:V\To\P^1\times\P^1\buildrel{\pr_1}\over\To\P^1$$ defined over $K$ with $V$ a smooth projective $K$-variety.
For a smooth projective curve $C$ and a ramified covering map $\gamma:C\To\P^1$, by pulling-back via $\gamma$ we obtain a Ch\^atelet surface bundles over $C$ fitting into the following commutative diagram of projective varieties
$$\xymatrix{
X\ar[r]\ar[d]\ar@/^2pc/[dd]^(.25){F}&V\ar[d]\ar@/^2pc/[dd]^(.25){f}    \\
C\times\P^1\ar[r]\ar[d]^{\pr}&\P^1\times\P^1\ar[d]^{\pr_1} \\
C\ar[r]^\gamma&\P^1
}$$
It follows immediately from Lemma \ref{smoothfiberlemma} that the projective variety $X$ is smooth once the branch locus $\mathfrak{B}_\gamma$  of $\gamma:C\To\P^1$ is disjoint from the degeneracy locus $\mathfrak{B}$ of $f:V\To \P^1$ (as defined in Definition \ref{degeneracylocusofChsurfbun}).

Moreover, under  mild conditions the Brauer group of $X$ can be well-understood, which plays an important role in the study of the arithmetic of $X$.

\begin{prop}[Poonen et al.]\label{propBr}
Suppose that $\textup{De}(V)$ is a smooth curve and that
%the branch locus
$\mathfrak{B}_\gamma$
%of $\gamma:C\To\P^1$
is disjoint from $\mathfrak{B}$, then $$\Br(C)\To\Br(X)$$ is a surjection.
\end{prop}

\begin{proof}
By explicit calculations of the geometric Picard groups  \cite[Lemma 5.1]{Poonen}, Poonen proved the surjectivity of $\Br(C\times\P^1)\To\Br(X)$, from which the surjectivity in the statement is hence deduced.

In \cite[Proposition 2.1]{CTsurPoonen}, Colliot-Th\'el\`ene gave an alternative proof of the surjectivity mentioned above.
Notice that the degeneracy locus $\textup{De}(X)$ of the conic bundle $X\To C\times\P^1$ is the curve $C\times_{\P^1}\textup{De}(V)$. It is a smooth geometrically integral curve under our assumption, cf. \cite[Lemma 7.1]{Poonen}. Sitting above the generic point of $\textup{De}(X)$ is  the only possible non-split fiber over a 1-codimensional point of $C\times\P^1$. A study of the purity exact sequence of the Brauer group allows us to conclude. See also \cite[Proposition 2.2]{CtPalSk}.
\end{proof}

We also know that all \'etale covers of $X$ arise from $C$.
\begin{lem}[{\cite[Lemma 8.1]{Poonen}}]\label{etalecovers}
The morphism $F:X\To C$ induces an equivalence of categories of finite \'etale covers $\textup{\textbf{Et}}(C)\To\textup{\textbf{Et}}(X)$.
\end{lem}

In \cite{Poonen}, Poonen chose an admissible polynomial $P_t(x)$ with degree $d=2$ in the variable $t$ such that the Ch\^atelet surface fiber $V_\infty$ violates Hasse principle. For any curve $C$ possessing only finitely many $K$-rational points, there exists  a covering map $\gamma$ sending all such $K$-points to $\infty\in\P^1$ and such that the branch locus $\mathfrak{B}_\gamma$ of $\gamma$ is disjoint from $\mathfrak{B}$. Now it is clear that the Ch\^atelet surface bundle  $X$ violates the Hasse principle. Moreover Proposition \ref{propBr} allows him to show that this violation of Hasse principle cannot be explained by the Brauer\textendash Manin obstruction (even applied to \'etale covers).

\subsection{Lagrange interpolation applied for algebraic fibrations}

We apply Lagrange interpolation to construct  Ch\^atelet surface bundles over $\mathbb{P}^1$. This is the key ingredient of the proof of our main theorem.

\begin{prop}[Interpolation of Ch\^atelet surfaces]\label{interpolation-GeneralChSur}
Given a finite set of closed points $\{\theta_i\}_{i\in I}$ of $\mathbb{A}^1$ of residue field $K(\theta_i)$ and given for each $i\in I$ a Ch\^atelet surface $\mathcal{S}_i$ defined over $K(\theta_i)$, then there exists a Ch\^atelet surface bundle over $\mathbb{A}^1$ defined over $K$ whose fiber at $\theta_i$ is isomorphic to $\mathcal{S}_i$.
\end{prop}

\begin{proof}
Let $\mathcal{S}_i$ be defined by $y^2-a_iz^2=P_i(x)$ where $a_i\in K(\theta_i)$ is a nonzero constant and  $P_i(x)\in K(\theta_i)[x]$ is a separable polynomial of degree $4$. We are looking for a polynomial $a_t\in K[t]$ and a polynomial $P_t(x)\in K[t,x]$ such that for all $i\in I$ we have the evaluations $a_{\theta_i}=a_i$ and $P_{\theta_i}(x)=P_i(x)$. Then the Ch\^atelet bundle surface bundle over $\mathbb{A}^1$ defined by $$y^2-a_tz^2=P_t(x)$$  has $\mathcal{S}_i$ as its fiber at $\theta_i$.

If we write $P_t(x)$ as $A_tx^4+B_tx^3+C_tx^2+D_tx+E_t$ with $A_t,B_t,C_t,D_t,E_t\in K[t]$ and write $P_i(x)$ as $A_ix^4+B_ix^3+C_ix^2+D_ix+E_i$ with $A_i,B_i,C_i,D_i,E_i\in K(\theta_i)$. Our aim reduces to a question of interpolating given values at the $K(\theta_i)$-rational points $\theta_i$.

Each closed point $\theta_i$ corresponds to a unique monic irreducible polynomial $\varphi_i\in K[t]$ of degree $d_i$. The residue field $K(\theta_i)=K[t]/(\varphi_i)$. These $\varphi_i$'s are pairwise prime to each other. Consider the evaluation morphism
$$\textup{ev}:K[t]\To\prod_{i\in I}K(\theta_i)$$ mapping $g\in K[t]$ to $(g(\theta_i))_{i\in I}=(g\mod (\varphi_i))_{i\in I}$.
By Chinese remainder theorem, the morphism $\textup{ev}$ is surjective and its kernel is the ideal of $K[t]$ generated by the polynomial $\varphi=\prod_{i\in I}\varphi_i$. Whence, for any given family of algebraic numbers $a_i$, $A_i$, $B_i$, $C_i$, $D_i$, and $E_i$ in $K(\theta_i)$ the desired polynomials $a_t$, $A_t$, $B_t$, $C_t$, $D_t$, and $E_t$ exist and they are unique modulo $\varphi\in K[t]$.
\end{proof}

As discussed in \S \ref{ChSurBun-section}, if we want to extend the Ch\^atelet surface bundle  $V^0\to\mathbb{A}^1$ obtained above to a bundle $V$ over $\P^1$, one sufficient condition is that $a_t$ is a constant and that $P_t(x)=A_tx^4+B_tx^3+C_tx^2+D_tx+E_t$ is of even degree $d$ in the variable $t$. Furthermore, if we want to apply Proposition \ref{propBr} to study the arithmetic of the bundle obtained by pulling-back $V$, we would like to choose an admissible $P_t(x)$ so that $\textup{De}(V)$ is a smooth curve. However, it seems not easy to discuss the admissibility of $P_t(x)$ in such a general setting where the given data $\{\mathcal{S}_i\}_{i\in I}$    are allowed to be arbitrary Ch\^atelet surfaces. It turns out that, for our purpose, it suffices to make the following further restrictions.

\begin{prop}\label{interpolation-SpecialChSur}
Given a finite set of closed points $\{\theta_i\}_{i\in I}$ of $\mathbb{P}^1$ of residue field $K(\theta_i)$ and given for each $i\in I$ a Ch\^atelet surface $\mathcal{S}_i$ over $K(\theta_i)$ defined by the equation
$$y^2-a_iz^2=A_ix^4+B_ix^3+C_ix^2+D_ix+E_i$$ with $a_i,A_i,B_i,C_i,D_i,E_i\in K(\theta_i)$.
Suppose further that for all $i\in I$
\begin{itemize}
\item[-] the nonzero constants $a_i=a\in K$ are the same and they lie in $K$,
\item[-] the constants $B_i=D_i=0$,
\item[-] the constant $E_i\neq0$.
\end{itemize}

Then there exists a Ch\^atelet surface bundle $V\To \P^1$ interpolating $\mathcal{S}_i$, i.e. $V_{\theta_i}\simeq\mathcal{S}_i$ for all $i\in I$, such that both $V$ and $\textup{De}(V)$ are smooth over $K$.
\end{prop}

\begin{proof}
It is clear that the set of points $\{\theta_i\}_{i\in I}$ is contained in a certain open subset of $\P^1$ which is isomorphic to $\mathbb{A}^1$.
According to Proposition \ref{interpolation-GeneralChSur} and its proof, we obtain a polynomial $P_t(x)\in K[t,x]$ giving a Ch\^atelet surface bundle $V^0\to\mathbb{A}^1$ interpolating $\mathcal{S}_i$. Our strategy is to modify $P_t(x)$ without changing the fibers of $V^0$ over $\theta_i$. We need to find an admissible $P_t(x)$ of even degree $d$ in the variable $t$ so that $V^0\to\mathbb{A}^1$ can be extended to a Ch\^atelet surface bundle $V\to\P^1$ such that $\textup{De}(V)$ is smooth. Then $V$ is also smooth over $K$ by Lemma \ref{totalsmoothnesslemma}.

Under the assumptions, the polynomial $a_t=a$ is a nonzero constant in $K$ and the polynomial $P_t(x)$ is of the form $A_tx^4+C_tx^2+E_t$.
Let $\varphi$ be as in the proof of Proposition \ref{interpolation-GeneralChSur}.
We observe that if we replace $A_t$ by $A_t+\varphi\psi^A$ with an arbitrary $\psi^A\in K[t]$ we will obtain another Ch\^atelet surface bundle over $\mathbb{A}^1$ without changing the fiber over each $\theta_i$. The same observation applies to $C_t$ and $E_t$ as well. We will prove the existence of suitable $\psi^A,\psi^C,\psi^E\in K[t]$, i.e. a suitable modification of the bundle, such that the further requirement of the admissibility of $P_t(x)$ will be satisfied.

First of all, if the degree of $\psi^A$ is large enough, then $$\deg(A_t+\varphi\psi^A)=\deg(\varphi)+\deg(\psi^A).$$ Without lost of generality, we may assume that $$\deg{\psi^A}=\deg{\psi^C}=\deg{\psi^E}=d_0$$ and that $$d_0+\deg{\varphi}=d$$ is even. From now on, we fix the degree and we still have enough freedom to choose suitable $\psi^A,\psi^C,\psi^E\in K[t]$. After the modifications, we will have $$\deg(A_t)=\deg(C_t)=\deg(E_t)=d.$$
According to the discussion in \S \ref{ChSurBun-section},  we are able to extend the modified bundle to a Ch\^atelet surface bundle over $\P^1$.

Secondly, we are going to find out some sufficient conditions on the modified $A_t$, $C_t$ and $E_t$ that will imply the admissibility in the  statement. We denote by $\Delta_t=C^2_t-4A_tE_t$ the ``discriminant'' of $P_t(x)$.
%For any algebraic number $\theta\in\bar{K}$ such that $A_\theta\neq0$, the quartic polynomial $P_\theta(x)$ is separable if and only if $\Delta_\theta\cdot E_\theta\neq0$. Whence, both $P_\theta(x)$ and $P^*_\theta(x')=x'^4P_\theta(1/x')$ are separable if and only if $A_\theta\cdot\Delta_\theta\cdot E_\theta\neq0$. (not useful??)

Recall by definition that $P_t(x)$ is admissible if $\textup{De}(V)\subset\P^1\times\P^1$ defined by $$\tilde{P}_{t,t'}(x,x')=0$$ is a smooth curve.
By definition
\begin{equation*}
    \begin{array}{ll}
        \tilde{P}_{t,t'}(x,x') & =t'^dx'^4P_{t/t'}(x/x')\\
                  & =t'^dA_{t/t'}x^4+t'^dC_{t/t'}x^2x'^2+t'^dE_{t/t'}x'^4.
    \end{array}
\end{equation*}
Also recall that $\P^1\times\P^1$ is covered by four affine open subsets $\mathcal{U}_{\infty\infty}$, $\mathcal{U}_{\infty0}$, $\mathcal{U}_{0\infty}$, and $\mathcal{U}_{00}$, cf. \S \ref{ChSurBun-section} for the notation. To check the smoothness of $\textup{De}(V)$, we need to restrict to each of these four affine plans.
Since the form of $\tilde{P}_{t,t'}(x,x')$ is symmetric between $t$ and $t'$, and between $x$ and $x'$, we only need to study its restriction to $\mathcal{U}=\mathcal{U}_{\infty\infty}$ and can easily find out the right condition for the other three.

The affine curve $\textup{De}(V)\cap \mathcal{U}\subset\mathbb{A}^2$ is given by $$P_t(x)=A_tx^4+C_tx^2+E_t=0.$$
We fix an algebraic closure $\bar{K}$ of $K$. By Jacobian criterion, the coordinates (also denoted by $(t,x)\in\bar{K}^2$ by abuse of notation) of its singular geometric points satisfy simultaneously the following equations.
\begin{equation*}
    \begin{array}{ll}
        (1)& 0 =  A_tx^4+C_tx^2+E_t\\
        (2)&0 =  A'_tx^4+C'_tx^2+E'_t\\
        (3)&0 = 4A_tx^3+2B_tx
    \end{array}
\end{equation*}
Three cases may happen.
\begin{itemize}
\item[-] Suppose that $x=0$. Then we must have simultaneously $E_t=0$ and $E'_t=0$, i.e. the polynomial $E_t$ has double roots in $\bar{K}$.
\item[-] Suppose that $x\neq0$ and $A_t=0$. Then we must have $C_t=0$ and $E_t=0$, i.e. the polynomial $A_t$, $C_t$ and $E_t$ have at least one common root in $\bar{K}$.
\item[-] Suppose that $x\neq0$ and $A_t\neq0$. We set $X=x^2$, the equation (3) implies that $X=x^2=-C_t/2A_t$. Then $A_tX^2+C_tX+E_t=0$ must have a double root in $X$. Therefore $\Delta_t=0$ and $x^4=E_t/A_t$. Substitute $x^2$ and $x^4$ to equation (2), we obtain $\Delta'_t=0$. Hence $\Delta_t$ has double roots in $\bar{K}$.
\end{itemize}
In other words, if the following three conditions are satisfied, the affine curve $\textup{De}(V)\cap \mathcal{U}$ is smooth.
\begin{itemize}
\item[-] \quad $A_t$ is separable,
\item[-] \quad $\Delta_t=C^2_t-4A_tE_t$ is separable,
\item[-] \quad $A_t$, $C_t$, $E_t$ have no common roots in $\bar{K}$.
\end{itemize}
Concerning the symmetry between $x$ and $x'$, we should add the separability of $E_t$ to the first condition and the last two conditions stay unchanged.
Concerning the symmetry between $t$ and $t'$, we need to consider reciprocal polynomials of $A_t$, $C_t$, $E_t$, and $\Delta_t$. It is clear that the separability of $A_t$ implies the separability of $A^*_{t'}$. One can check that the reciprocal polynomial $\Delta^*_{t'}$ of $\Delta_t$ (of degree $2d$) is equal to $C^{*2}_{t'}-4A^*_{t'}E^*_{t'}$. So its separability yields no new condition. If $A_t$, $C_t$, $E_t$ have no common roots, then $A^*_{t'}$, $C^*_{t'}$, $E^*_{t'}$ have no common roots either.

Hence, if the following three conditions are simultaneously satisfied, then $\textup{De}(V)\subset\P^1\times\P^1$ is smooth and $P_t(x)$ is admissible.
\begin{enumerate}
\item $A_t$ and $E_t$ are separable,
\item $\Delta_t=C^2_t-4A_tE_t$ is separable,
\item no irreducible polynomial in $K[t]$ can divide $A_t$, $C_t$, and $E_t$ at the same time.
\end{enumerate}

Finally, we are going to prove the existence of $\psi^A$, $\psi^C$ and $\psi^E$ such that the modified $A_t$, $C_t$ and $E_t$ satisfy the conditions above.

\textbf{Choice of $C_t$}

The modified $A_t$, $C_t$ and $E_t$ will be $A_t+\varphi\psi^A$, $C_t+\varphi\psi^C$ and $E_t+\varphi\psi^E$ with $\psi^A,\psi^C,\psi^E\in K[t]$ of the same large enough degree. The polynomial $\psi^C$ can be chosen arbitrarily, and from now on $\psi^C$ and $C_t^{\textup{new}}=C_t+\varphi\psi^C$ are fixed.

\textbf{Choice of $E_t$}

For $n\in\mathbb{N}$, consider $${E}_t^{(n)}=\frac{1}{n}E_t+\varphi\psi_0^E$$ with $\psi_0^E\in K[t]$ separable and prime to $\varphi\cdot C^\textup{new}_t$.  When $n$ is large, the polynomial $ {E}^{(n)}_t$ is close enough to the separable polynomial $\varphi\psi_0^E$ with respect to the topology induced by a certain embedding $K\To \mathbb{C}$, it must also be  separable by the forthcoming Lemma \ref{separablelemma}.

By the construction of $E_t$ in  Proposition \ref{interpolation-GeneralChSur}, we know that $E_t\mod\varphi_i=E_i\in K(\theta_i)$ for all $i\in I$. The assumption that the constant $E_i\neq0$ implies that the irreducible polynomial $\varphi_i$ does not divide $E_t$. Therefore $\varphi$ is prime to $E_t$. Let $q_t\in K[t]$ be an irreducible polynomial dividing $C^\textup{new}_t$.
For any positive integers $k<k'$, if $q_t$ divides both $ {E}^{(n+k)}_t$ and $ {E}^{(n+k')}_t$, then it divides $$\varphi\psi_0^E=\frac{1}{k'-k}[(n+k') {E}^{(n+k')}_t- (n+k){E}^{(n+k)}_t].$$ According to our choice of $\psi_0^E$, it is prime to $C^\textup{new}_t$. So $q_t$ divides $\varphi$. But then it divides $$E_t=(n+k)( {E}^{(n+k)}_t-\varphi\psi_0^E)$$ contradicting the fact that $\varphi$ and $E_t$ are prime to each other. Hence $C_t$ together with any two polynomials in the set $\Sigma=\{ {E}_t^{(n+k)}|k\in\mathbb{N}\}$  cannot have a common irreducible divisor. Now $C_t$ has only finitely many irreducible divisors but $\Sigma$ is an infinite set, by the pigeonhole principle there exists $k\in\mathbb{N}$ such that $C_t$ and $ {E}_t^{(n+k)}$ are prime to each other.

Now we take $$\psi^E=(n+k)\psi_0^E$$ and
\begin{equation*}
    \begin{array}{ll}
        E_t^\textup{new} & =(n+k){E}^{(n+k)}_t\\
                  & =E_t+(n+k)\varphi\psi_0^E\\
        & =E_t+\varphi\psi^E,
            \end{array}
\end{equation*}
then $E_t^\textup{new}$ is separable and prime to $C_t$. The condition (3) mentioned above is already satisfied for the chosen $C_t$, $E_t$ and arbitrary $A_t$.

\textbf{Choice of $A_t$}

For $n\in\mathbb{N}$, consider $${A}^{(n)}_t=\frac{1}{n}A_t+\varphi\psi_0^A$$ with $\psi_0^A$ prime to $\varphi\cdot E^\textup{new}_t$.
According to Lemma \ref{separablelemma}, we will set $\psi^A=n\psi_0^A$ as well as
\begin{equation*}
    \begin{array}{ll}
        A^\textup{new}_t & =n{A}^{(n)}_t\\
                  & =A_t+n\varphi\psi_0^A\\
        & =A_t+\varphi\psi^A
    \end{array}
\end{equation*}
which is a separable polynomial if $n$ is chosen large enough. Condition (1) is satisfied.

By our construction $E^\textup{new}_t\equiv E_t\mod\varphi$. The assumption that $E_i\neq0$ for all $i\in I$ implies that $\varphi$ and $E^\textup{new}_t$ are prime to each other.

We have
\begin{equation*}
    \begin{array}{ll}
        \Delta^\textup{new}_t& =C^{\textup{new}2}_t-4A^\textup{new}_tE^{\textup{new}}_t\\
                  & =C^{\textup{new}2}_t-4(A_t+n\varphi\psi_0^A)E^{\textup{new}}_t\\
        & =(C^{\textup{new}2}_t-4A_tE^{\textup{new}}_t)-4n\varphi\psi_0^AE^\textup{new}_t.
    \end{array}
\end{equation*}
Since $\varphi\psi_0^AE^\textup{new}_t$ is a separable polynomial,  Lemma \ref{separablelemma} implies that $\Delta^\textup{new}_t$ is also separable if $n$ is large enough. Condition (2) is satisfied.

In summary, the polynomial $$P^\textup{new}_t(x)=A^\textup{new}_tx^4+C^\textup{new}_tx^2+E^\textup{new}_t\in K[t,x]$$ is admissible of even degree in $t$. The equation $y^2-az^2=P^\textup{new}_t(x)$ defines a Ch\^atelet surface bundle $V\To\P^1$ over $K$ such that for each $i\in I$ the fiber $V_{\theta_i}$ is the given Ch\^atelet surface $\mathcal{S}_i$.
\end{proof}

\begin{lem}\label{separablelemma}
Let $f\in\C[t]$ be a separable polynomial of degree $d>1$. Consider the product topology on the subset $\C_d[t]$ of polynomials of degree at most $d$. If $g\in\C_d[t]$ is close enough to $f$, then $g$ is also separable.

In particular, for any $f_0\in\C_d[t]$ the polynomial $f_0+nf$ is separable for $n\in\mathbb{N}$ large enough.
\end{lem}

\begin{proof}
Let $\Delta(h)$ be the discriminant of a degree $d$ polynomial $h=a_dt^d+\cdots+a_2t^2+a_1t+a_0$. It is a polynomial in $n+1$ variables $a_d,\ldots,a_2,a_1,a_0$, it is a continuous function under the complex topology. $\Delta(h)$ is nonzero if and only if $h$ is separable.
Since $\Delta(f)\neq0$ and $g$ is close enough to $f$, we find that $\Delta(g)\neq0$ and therefore $g$ is also separable.

In particular, for $n\in\mathbb{N}$ such that $g=\frac{1}{n}f_0+f$ is close enough to $f$, it is separable. So is $ng=f_0+nf$.
\end{proof}

\subsection{Proof of the main theorem}

We are ready to prove the main theorem \ref{maintheorem}.
\subsubsection{Construction of the curve $C$ and $\gamma:C\To\P^1$.}\ \\

Let $K$ be a number field.
Let $E\subset\P^2$ be an elliptic curve given by the  Weierstrass equation
$$y^2=x^3+ax+b$$
where $a,b\in K$ with $4a^3+27b^2\neq0$.
For a non-constant non-linear polynomial $h(t)\in K[t]$, we put $H(t)=h(t)+ah(t)^3+bh(t)^4\in K[t]$. With the convention in \S \ref{hyperelliptic}, the equation $$s^2=H(t)=h(t)+ah(t)^3+bh(t)^4$$ defines a smooth projective hyperelliptic curve $C$ over $K$ of genus $g\geq2$ once $H(t)$ is a separable polynomial. It is  a double cover $\gamma:C\to\P^1$ of $\P^1$ via the projection to $t$, whose branch locus $\mathfrak{B}_\gamma\subset \P^1$  consists of  finitely many closed points.

\begin{lem}\label{CtoE}
Assume that the curve $C$ is smooth, for example when $H(t)$ is separable.
Then there exists a non-constant morphism of algebraic curves $$\phi:C\To E.$$
\end{lem}

\begin{proof}
Consider the open subset $C^0$ of $C$ given by the affine equation in $\mathbb{A}^2$
$$s^2=H(t)=\left(1+ah(t)^2+bh(t)^3\right)h(t)$$
with $1+ah(t)^2+h(t)^3\neq0$.
Then we can define a morphism $$\phi^0:C^0\To \P^2$$ by $$(s,t)\mapsto\left(s:1+ah(t)^2+h(t)^3:h(t)s\right).$$
Since $$s^2=\left(1+ah(t)^2+bh(t)^3\right)h(t)$$ on $C^0$, we have
$$\left(1+ah(t)^2+bh(t)^3\right)^2h(t)s=\left(1+ah(t)^2+bh(t)^3\right)s^3$$
or equivalently
$$\left(1+ah(t)^2+bh(t)^3\right)^2\left(h(t)s\right)=s^3+as\left(h(t)s\right)^2+b\left(h(t)s\right)^3,$$
which signifies that the image of $\phi^0$ lies on the elliptic curve $E$ defined by the homogeneous Weierstrass equation $$y^2z=x^3+axz^2+bz^3.$$
By the valuative criterion of properness, the morphism $\phi^0$ extends uniquely to a non-constant morphism $\phi:C\To E$ between smooth projective curves.
\end{proof}

\begin{prop}\label{curveC}
Let $L/K$ be a nontrivial extension of number fields.
Assume that there exists an elliptic curve $E$ defined over $K$  of Mordell\textendash Weil rank $0$ such that the subgroup of divisible elements of its  Tate\textendash Shafarevich group is trivial.

Then there exists a hyperelliptic curve $C$ of genus at least $2$ such that
\begin{itemize}
\item[-] $C(K)\neq\varnothing$,
\item[-] there exists a closed point $\Theta'$ of $C$ whose image $\Theta=\gamma(\Theta')\in\P^1$ is of residue field $L$, in particular $C(K)\subsetneq C(L)$,
\item[-] for any non-archimedean place $v_0\in\Omega_K$ and for any $(x_v)_{v\in\Omega_K}\in C(\textup{\textbf{A}}_K)^\Br$, the $v_0$ component $x_{v_0}$ lies inside the image of $C(K)\to C(K_{v_0})$.
\end{itemize}
\end{prop}

\begin{proof}
Let $E$ be the elliptic curve defined by $y^2=x^3+ax+b$ satisfying the hypotheses mentioned in the statement.
Under the assumption that $H(t)$ is separable, let $C$ be the hyperelliptic curve defined by
$$s^2=H(t)=h(t)+ah(t)^3+bh(t)^4.$$
Apply \cite[Theorem 8.6]{Stoll07} to the non-constant morphism  $\phi:C\to E$ obtained in Lemma \ref{CtoE}, we find that $C$ is  excellent with respect to abelian coverings in the sense of M. Stoll, cf. \cite[Definition 6.1(6)]{Stoll07}.
Then the last assertion follows from \cite[Corollary 7.3]{Stoll07}.

If $h(t)\in K[t]$ has a linear factor $t-c_0$, then the point of coordinates $(0,c_0)$ is a $K$-rational point of $C$. If $h(t)$ has an irreducible factor $\varphi(t)$ such that $L=K[t]/(\varphi)$, then the closed point $\Theta'\in C$ of coordinates $(0,\vartheta)$ is an $L$-rational point where $\vartheta$ denotes the class of $t$. And its image $\Theta=\gamma(\Theta')\in\P^1$ admitting coordinate $\vartheta$ is of residue field $K(\vartheta)=L$.
According to the formula of $\phi:C \to E$ given in the proof of Lemma \ref{CtoE}, we remark that these points of $C$ lie above  the identity element $(0:1:0)$ of the elliptic curve $E$.

We fix an irreducible polynomial $\varphi$ such that $L=K[t]/(\varphi)$. We could set $h(t)=t\varphi(t)$, but then $H(t)$ is not separable in general.
It remains to modify $h(t)$ to ensure the separability of $H(t)$ and thus the smoothness of $C$.
Since $L$ is a non-trivial extension of $K$, the polynomial $\varphi$ has no roots in $K$. Instead of $h(t)$, consider a family of separable polynomials $$h_c(t)=\frac{(t+c)\varphi(t+c)}{c\varphi(c)}=\frac{h(t+c)}{h(c)}$$ parameterized by $c\in K^*$. It is clear that each member of the family has a factorisation as required in the previous paragraph, thus it suffices to choose a constant $c\in K^*$ such that $H_c(t)=h_c(t)\left(bh_c(t)^3+ah_c(t)^2+1\right)$ is separable.

When $b\neq0$, the cubic polynomial $bx^3+ax^2+1$ has $-4a^3-27b^2\neq0$ as its discriminant, thus  factorises as $$bx^3+ax^2+1=b(x-\alpha_1)(x-\alpha_2)(x-\alpha_3)$$
with distinct roots $\alpha_1,\alpha_2,\alpha_3\in \bar{K}^*$ where $\bar{K}$ denotes an algebraic closure of $K$.
Now $$H_c(t)=bh_c(t)\left(h_c(t)-\alpha_1\right)\left(h_c(t)-\alpha_2\right)\left(h_c(t)-\alpha_3\right)$$ is separable if  $h_c(t)-\alpha_i$ is separable for each $i\in\{1,2,3\}$.

We conclude by claiming that for any $\alpha\in\bar{K}^*$ there exists at most finitely many $c\in K^*$ such that $h_c(t)-\alpha$ is not separable. Indeed, the separability of $h_c(t)-\alpha$ is equivalent to that of $h(t+c)-\alpha h(c)$.
By taking derivative, we see that $h(t+c)-\alpha h(c)$ is not separable  only if $$\alpha h(c)\in \{h(\lambda);~h'(\lambda)=0,~\lambda\in\bar{K}\},$$ which is possible for at most finitely many $c\in K^*$.

For the case where $b=0$ we find that $a\neq0$. The same argument  applied to $ax^2+1$ allows to conclude as well.
\end{proof}

\begin{rem} According to \cite[Theorem 1.1]{MazurRubin10}, there exists always an elliptic curve $E$ defined over $K$ such that $E(K)$ is finite, the relevant assumption is the triviality of the subgroup of divisible elements of the Tate\textendash Shafarevich group of $E$.
When $K=\Q$ many elliptic curves are proved to have finite Tate\textendash Shafarevich group and of Mordell\textendash Weil rank $0$, for example the curve with LMFDB \cite{LMFDB} label 64.a3 defined by $$y^2=x^3-4x.$$
It can  be deduced that over some other fields $K$ the assumption of Proposition \ref{curveC} or Theorem \ref{maintheorem} can be satisfied. For example $K=\Q(\sqrt{-1})$, $\Q(\sqrt{3})$, and many others, please refer to \cite[Table 1 of \S4.6 and Lemma 4.7]{Liang-noninv}.
\end{rem}

\subsubsection{Construction of the bundle $V\To\P^1$.}\ \\

According to Faltings' theorem, the set $C(L)$ of $L$-rational points on $C$ is finite.
The subset $\gamma(C(L))\subset \P^1(L)$ corresponds to a finite set $\mathfrak{R}$ of closed points of $\P^1$.
It follows from Proposition \ref{curveC} that $\mathfrak{R}$ contains at least one point $\Theta=\gamma(\Theta')$ of residue field $L$ and at least one point of residue field $K$, where the latter is denoted by $O=\gamma(O')$ with a certain $O'\in C(K)$.
The set $\mathfrak{R}$ may also contain other closed points of residue field a certain intermediate field of $L/K$.

For each point  $\theta\in\mathfrak{R}\cup\mathfrak{B}_\gamma$, we will choose a Ch\^atalet surface $\mathcal{S}_\theta$ then apply the interpolation method of Proposition \ref{interpolation-SpecialChSur}.

Above all, we are going to fix a constant $a\in \mathcal{O}_K$ for all Ch\^atalet surfaces that will appear here.
According to the Chebotarev density theorem, we choose an odd non-archimedean place $v_0$ of $K$ which splits completely in $L$. We also fix a place $w_0$ of $L$ lying over $v_0$.
By weak approximation of the number field $K$, one can find a nonzero element $a\in K$ such that
\begin{itemize}
\item[-] the valuation $v_0(a)$ is odd,
\item[-] the constant $a$ is positive for all real places of $K$,
\item[-] the constant $a$ is a square in $K_v$ for all even places $v$.
\end{itemize}
Since these conditions are stable under multiplication by a square in $K^*$, we may assume that $a\in\mathcal{O}_K$. Because $a$ is a nonzero $2\R$ local square such that $v_0(a)$ and $w_0(a)$ are odd, Propositions \ref{WuHan1} and \ref{WuHan2} can be applied.

The  Ch\^atelet surface $\mathcal{S}_\theta$ will be  defined over the residue field $K(\theta)$ by the equation $$y^2-az^2=P_\theta(x),$$ where $P_\theta(x)\in K(\theta)[x]$ is always a separable polynomial of the form $A_\theta x^4+C_\theta x^2+E_\theta$ with $E_\theta\neq0$ given according to different $\theta\in\mathfrak{R}\cup\mathfrak{B}_\gamma$ in the following cases.
\begin{itemize}
\item[-] For the particular point $\theta=O\in \mathfrak{R}$, the Ch\^atelet surface $\mathcal{S}_O$ defined over $K$ is obtained by Proposition \ref{WuHan1}. Then $\mathcal{S}_O(K_{v_0})=\varnothing$ and $\mathcal{S}_O(K_{v})\neq\varnothing$ for all $v\neq v_0$.
\item[-] For the particular point $\theta=\Theta\in \mathfrak{R}$, the Ch\^atelet surface $\mathcal{S}_\Theta$ defined over $L$ is obtained by Proposition \ref{WuHan2}. It violates Hasse principle, i.e. $\mathcal{S}_\Theta(L)=\varnothing$ but $\mathcal{S}_\Theta(L_w)\neq\varnothing$ for all places $w$ of $L$.
\item[-] For any $\theta\in \mathfrak{R}\setminus\{\Theta,O\}$, if $\mathfrak{R}\setminus\{\Theta,O\}$ is nonempty, the Ch\^atelet surface $\mathcal{S}_\theta$ defined over a certain intermediate field $K(\theta)=M$ of $L/K$ is obtained by taking base change from $K$ to $M$ of the surface in the conclusion of Proposition \ref{WuHan1}. Then $\mathcal{S}_\theta(L)=\varnothing$ since $\mathcal{S}_\theta(L_{w_0})=\varnothing$ for any place $w_0$ lying over $v_0$.
\item[-] For any $\theta\in \mathfrak{B}_\gamma\setminus\mathfrak{R}$, if $\mathfrak{B}_\gamma\setminus\mathfrak{R}$ is nonempty, the polynomial $P_\theta$ can be any separable polynomial with coefficients in $K(\theta)$ of the form $A_\theta x^4+C_\theta x^2+E_\theta$ with $E_\theta\neq0$, then $\mathcal{S}_\theta$ is smooth.
\end{itemize}

Now we apply the interpolation method in Proposition \ref{interpolation-SpecialChSur} to obtain a Ch\^atelet surface bundle $V\To\P^1$ such that $V_\theta\simeq\mathcal{S}_\theta$ for all $\theta\in\mathfrak{R}\cup\mathfrak{B}_\gamma$ and such that both $V$ and $\textup{De}(V)$ are smooth over $K$.

\subsubsection{Construction of $X\To C$ and conclusion.}\ \\

For all $\theta\in\mathfrak{R}\cup\mathfrak{B}_\gamma$ the fibers $V_\theta$ are smooth, hence the branch locus $\mathfrak{B}$ of $\pr_1:\textup{De}(V)\subset\P^1\times\P^1\To\P^1$ is disjoint from $\mathfrak{B}_\gamma$ by Lemma \ref{smoothfiberlemma} and then the pull-back $X$ of $V$ by $\gamma:C\To\P^1$ is a smooth $K$-variety.

\textbf{We study the arithmetic of $X$ over $L$.}

For any $\theta'\in C(L)$, its image $\theta=\gamma(\theta')\in \mathfrak{R}$. The fiber $X_{\theta'}=(V_\theta)_L$ does not possess any $L$-rational point. We deduce that $X(L)=\varnothing$.

When $\theta'=\Theta'$, we have $X_{\Theta'}=V_\Theta\simeq \mathcal{S}_\Theta$ because both residue fields of $\Theta$ and $\Theta'$ are $L$.
Then $X_{\Theta'}(L_w)\neq\varnothing$ for all places of $L$. Any family of local points on the  fiber $X_{\Theta'}$ is orthogonal to $\Br(C_L)$ since it projects to a single $L$-rational point of $C$. Since $\Br(C_L)\To\Br(X_L)$ is surjective by Proposition \ref{propBr}, the Brauer\textendash Manin set $X_L({\textup{\textbf{A}}}_L)^{\Br(X_L)}$ contains the non-empty set $\{\Theta'\}\times X_{\Theta'}({\textup{\textbf{A}}}_L)=\{\Theta'\}\times \mathcal{S}_{\Theta}({\textup{\textbf{A}}}_L)$. The failure of Hasse principle on $X_L$ is not explained by the \BM obstruction.

Let $G$ be a finite \'etale group scheme over $L$, and let $\eta:Y\To X_L$ be a right torsor under $G$. According to Lemma \ref{etalecovers}, the torsor $\eta$ arises from a right torsor (by abuse of notation) $\eta:\mathcal{C}\To C_L$ under $G$, i.e. the following commutative diagram on the left is a pull-back.
$$\xymatrix@!=2pc{
Y\ar[d]^{\mathcal{F}}\ar[r]^G_{\eta} & X_{L}\ar[d]^{F_L} &&& **[l]Y^\tau\times_{\mathcal{C}^\tau}D\subset Y^\tau\ar[d]^{\mathcal{F}^\tau}\ar[r]^{G^\tau}_{\eta^\tau} & X_{L}\ar[d]^{F_L}\ar[r]&V_L\ar[d]^{f_L}\\
\mathcal{C}\ar[r]^G_{\eta} &C_L &&& **[l]D\subset\mathcal{C}^\tau\ar[r]^{G^\tau}_{\eta^\tau} &C_L\ar[r]_{\gamma_L}&\P^1_L
}$$
By the theory of torsors (see for example \cite[\S 2.2]{Skbook}), there exists a unique cohomology class $[\tau]\in\textup{H}^1(L,G)$ such that the twist $\mathcal{C}^{\tau}$ has an $L$-rational point $\Theta''$ mapping to $\Theta'\in C_L(L)$. Let $D$ be the connected component of $\mathcal{C}^\tau$ containing $\Theta''$. Since the branch loci of $\gamma_L$ and $\gamma_L\circ\eta^\tau$ are the same
because $\eta^\tau$ is \'etale, the argument in the previous paragraph applied to $\mathcal{F}^\tau:Y^\tau\times_{\mathcal{C^\tau}}D\To D$ implies that $$Y^\tau({\textup{\textbf{A}}}_L)^{\Br(Y^\tau)}\supset(Y^\tau\times_{\mathcal{C^\tau}}D)({\textup{\textbf{A}}}_L)^{\Br(Y^\tau\times_{\mathcal{C^\tau}}D)}\supset\{\Theta''\}\times (Y^\tau)_{\Theta''}({\textup{\textbf{A}}}_L).$$
Applying $\eta^\tau$ and taking union over all cohomology classes $[\sigma]\in\textup{H}^1(L,G),$ we obtain
$$\bigcup_{[\sigma]\in\textup{H}^1(L,G)}\eta^\sigma\left(Y^\sigma({\textup{\textbf{A}}}_L)^{\Br(Y^\sigma)}\right)\supset\eta^\tau\left(Y^\tau({\textup{\textbf{A}}}_L)^{\Br(Y^\tau)}\right)\supset\{\Theta'\}\times X_{\Theta'}({\textup{\textbf{A}}}_L).$$
Whence $$X_L({\textup{\textbf{A}}}_L)^{\textup{et},\Br}\supset \{\Theta'\}\times X_{\Theta'}({\textup{\textbf{A}}}_L)=\{\Theta'\}\times \mathcal{S}_{\Theta}({\textup{\textbf{A}}}_L)\neq\varnothing$$
by considering all finite \'etale group scheme $G$ over $L$ and all torsors under $G$. The failure of Hasse principle on $X_L$ is not explained by the \'etale \BM obstruction.

\textbf{We study the arithmetic of $X$ over $K$.}

As $X(L)=\varnothing$, a fortiori $X(K)=\varnothing$. For any place $v$ of $K$ different from $v_0$, it is clear that $X_{O'}=V_O\simeq\mathcal{S}_O$ possesses $K_v$-rational points, hence $X(K_v)\neq\varnothing$. Let $w_0$ be a place of $L$ lying over $v_0$. Since $v_0$ splits completely in $L$, the local fields $L_{w_0}$ and $K_{v_0}$ are equal. So $X(K_{v_0})=X(L_{w_0})\neq\varnothing$. The variety $X$ violates the Hasse principle. If there exists a family of local points $(x_v)_{v\in\Omega_K}$ orthogonal to $\Br(X)$, then $(F(x_v))_{v\in\Omega_K}$ is orthogonal to $\Br(C)$. We apply Proposition \ref{curveC} to find a $K$-rational point $\theta'$ of $C$ such that $\theta'=F(x_{v_0})$. However, the fiber $X_{\theta'}=V_\theta\simeq\mathcal{S}_\theta$ with $\theta=\gamma(\theta')\in \mathfrak{R}\cap\P^1(K)$ does not admit any $K_{v_0}$-rational point, which leads to a contradiction. So the \BM set of $X$ is empty,  the failure of the Hasse principle  over $K$ is accounted for by the \BM obstruction, a fortiori by the \'etale \BM obstruction as well.

The proof of Theorem \ref{maintheorem} is now completed.

\section{Further discussion}\label{further}

People naturally ask, for a nontrivial extension $L/K$ of number fields, whether there exists $X$  in a certain given class of varieties such that
\begin{itemize}
\item[-] \quad $\varnothing=X(K)\subsetneq X(\textup{\textbf{A}}_K)^{\Br(X)}\subset X(\textup{\textbf{A}}_K)$,
\item[-] \quad $\varnothing=X_L(L)= X_L(\textup{\textbf{A}}_L)^{\Br(X_L)}\subsetneq X_L(\textup{\textbf{A}}_L)$.
\end{itemize}
In other words, the violation of Hasse principle is not accounted for by the Brauer\textendash Manin obstruction  over $K$ but it becomes the case over $L$. The same question for the \'etale Brauer\textendash Manin obstruction also arises naturally. It is in some sense an opposite situation of our main Theorem \ref{maintheorem}.

Indeed, it never happens. For  smooth quasi-projective varieties $X$, we have
$$X(\textup{\textbf{A}}_K)^{\Br}\neq\varnothing \Longrightarrow X_L(\textup{\textbf{A}}_L)^{\Br}\neq\varnothing,$$
$$X(\textup{\textbf{A}}_K)^{\textup{et},\Br}\neq\varnothing \Longrightarrow X_L(\textup{\textbf{A}}_L)^{\textup{et},\Br}\neq\varnothing.$$

When $X$ is \emph{projective}, a sketch of proof of the first inclusion with the help of Weil restrictions has already appeared in a paper \cite[Remark 5, pages 95-96]{CTPoonen00} of J.-L. Colliot-Th\'el\`ene and B. Poonen. A proof with details is written in \cite[Lemma 2.1(1)]{CreutzViray21} by B. Creutz and B. Viray. Actually, they show that $X(\A_K)^{\Br(X)}\subset X_L(\A_L)^{\Br(X_L)}$. But this proof can not be extended directly to a proof for the second implication. Though we can apply the first implication to  \'etale covers of $X_L$, but the problem that not all $X_L$-torsors under a finite \'etale $L$-group scheme arise from a similar setting over $K$ prevents us from arriving at the conclusion.  Recently, Yang Cao has an argument that also makes use of Weil restrictions, but in another way such that it applies to both inclusions and such that it works more generally for quasi-projective varieties, cf. \cite{Cao2022}.

We should also mention that, for projective rationally connected varieties, there is one more alternative argument proposed to the second author by Colliot-Th\'el\`ene  using the theory of universal torsors.
If the Brauer\textendash Manin set over $K$ is non-empty, then there exist universal torsors of $X$, one of which possesses rational points everywhere locally, cf. \cite{CTSansuc87-1}. This torsor base extended to $L$ is a universal torsor of $X_L$ and it has rational points everywhere locally over $L$. Whence the Brauer\textendash Manin set over $L$ is not empty as well.

\bigskip

\begin{footnotesize}
%\noindent\textbf{Acknowledgements.}
% He would also like to thank  the referees for  suggestions on English expressions.
%The second author is partially supported by Anhui Initiative in Quantum Information Technologies No.AHY150200. Both authors are partially supported by the National Natural Science Foundation of China No.12071448.
\end{footnotesize}

\tableofcontents

\bigskip

%\nocite{*} %%show all the bib
\bibliographystyle{alpha}
\bibliography{mybib1}

\end{document}